\documentclass[a4paper,numberrefs]{icarnws}

\usepackage[latin2]{inputenc}
\usepackage[T1]{fontenc}
\usepackage{amsmath}
\usepackage{amssymb}
\usepackage{graphicx}
\usepackage{xspace}
\usepackage[mathscr]{eucal}
\usepackage{url}
\usepackage{subfig}
\usepackage{times}
\usepackage{tabularx}
\usepackage{amsthm}

\newcommand{\R}{\ensuremath{\mathbb{R}}}
\newcommand{\C}{\ensuremath{\mathbb{C}}}
\newcommand{\K}{\ensuremath{\mathbb{K}}}
\newcommand{\RICA}{\mbox{$\mathbb{R}$-ICA}\xspace}
\newcommand{\CICA}{\mbox{$\mathbb{C}$-ICA}\xspace}
\newcommand{\RISA}{\mbox{$\mathbb{R}$-ISA}\xspace}
\newcommand{\CISA}{\mbox{$\mathbb{C}$-ISA}\xspace}
\newcommand{\KICA}{\mbox{$\mathbb{K}$-ICA}\xspace}
\newcommand{\KISA}{\mbox{$\mathbb{K}$-ISA}\xspace}
\newcommand{\pv}{\ensuremath{\varphi}}
\newcommand{\G}{\ensuremath{\mathscr{G}}}
\newcommand{\F}{\ensuremath{\mathscr{F}}}

\renewcommand{\vec}{\mathbf}

\newtheorem{theorem}{Theorem}
\newtheorem*{proof2}{Proof}
\newtheorem{note}{Note}

\title{REAL AND COMPLEX INDEPENDENT SUBSPACE ANALYSIS BY GENERALIZED VARIANCE}

\author{
  \bf Zolt{\'a}n Szab{\'o} \hspace*{0.5cm} Andr{\'a}s L{\H{o}}rincz\\ 
  Department of Information Systems, E\"{o}tv\"{o}s Lor{\'a}nd University, \\
  P{\'a}zm{\'a}ny P. s{\'e}t{\'a}ny 1/C, Budapest H-1117, Hungary \\
  Research Group on Intelligent Information Systems \\
  Hungarian Academy of Sciences\\
  WWW home page: \url{http://nipg.inf.elte.hu}\\
  \tt szzoli@cs.elte.hu, andras.lorincz@elte.hu
}

\begin{document}
\maketitle

\begin{abstract}
Here, we address the problem of Independent Subspace Analysis (ISA). We develop a technique that (i) builds upon joint
decorrelation for a set of functions, (ii) can be related to kernel based techniques, (iii) can be interpreted as a
self-adjusting, self-grouping neural network solution, (iv) can be used both for real and for complex problems, and (v)
can be a first step towards large scale problems. Our numerical examples extend to a few 100 dimensional ISA tasks.
\end{abstract}

\keywords{Independent Subspace Analysis, joint \mbox{f-decorrelation}} 

\section{INTRODUCTION}
Uncovering independent processes is of high importance, because it breaks combinatorial explosion \citep{poczos06noncombinatorial}. In
cases, like Smart Dust, the problem is vital, because (i) elements have limited computational capacity and (ii) communication to remote
distances is prohibitively expensive. Self-adjusting, self-grouping neural network solutions may come to our rescue here. Here, we present
such an approach for Independent Subspace Analysis (ISA). The extension of ISA to Independent Process Analysis is straightforward under
certain conditions \citep{poczos06noncombinatorial}.

Our paper is built as follows. The \KISA model is introduced in Section~\ref{sec:KISA-model}.
Section~\ref{sec:KISA-method} is about our method. Illustrations are provided in Section~\ref{sec:illustrations}.

\section{THE \KISA MODEL}\label{sec:KISA-model}
Section~\ref{subsec:KISA-eqs} defines the \KISA task to be studied, Section~\ref{subsec:KISA-ambiguitites} treats the ambiguities of the
model.

\subsection{The \KISA Equations}\label{subsec:KISA-eqs}
We treat real and complex ISA tasks: Let $\K\in\{\R,\C\}$. Assume that we observe the mixture of multidimensional independent i.i.d.
sampled sources (\emph{components}):
\begin{align}
    \vec{z}(t)&=\vec{A}\vec{s}(t),&
    \vec{s}(t)&=[\vec{s}^1(t);\ldots;\vec{s}^M(t)],
\end{align} where $D=\sum_{m=1}^M d_m$ is the total dimension of the components, $\vec{A}\in\K^{D\times D}$ is the invertible \emph{mixing matrix}.
The task is to recover the hidden components \mbox{$\vec{s}^m(t)\in\K^{d_m}$} by means of observations  $\vec{z}(t)\in
\K^D$. If $\K=\R$ ($\K=\C$) then we shall talk about Real (Complex) ISA [i.e., \RISA (\CISA)] task. The special
$d_m=1(\forall m)$ case is the Real (Complex) Independent Component Analysis [i.e., \RICA (\CICA)].

\subsection{Ambiguities of the \KISA Model}\label{subsec:KISA-ambiguitites}
Identification of the $\KISA$ model is ambiguous. However, ambiguities are simple: hidden components $\vec{s}^m$ can be determined up to
permutation among subspaces and up to invertible transformation within subspaces. Details about \RISA and \CISA can be found in
\citep{theis04uniqueness1} and \citep{szabo06separation}, respectively.

Ambiguities within subspaces can be lessened: given our assumption on the invertibility of matrix $\vec{A}$, we can
assume without any loss of generality that both the sources and the observation are \emph{white}, that is,
\begin{eqnarray}
E[\vec{s}]&=&\vec{0},cov\left[\vec{s}\right]=\vec{I}_D,\label{eq:white1}\\
E[\vec{z}]&=&\vec{0},cov\left[\vec{z}\right]=\vec{I}_D,\label{eq:white2}
\end{eqnarray}
where $E[\cdot]$ denotes the expectation value, $\vec{I}_D$ is the \mbox{$D$-dimensional} identity matrix. Now, the
$\vec{s}^m$ sources are determined up to (i) permutation \emph{and} orthogonal transformation in the real case and (ii)
permutation \emph{and} unitary transformation in the complex case.

\section{\KISA BY JOINT DECORRELATION}\label{sec:KISA-method}
Components $\vec{s}^m$ are estimated by a neural network, which aims to `decorrelate' (see below) the
$\vec{y}^m\in\K^{d_m}$ parts of the $\K^D\ni\vec{y}(t)=[\vec{y}^1(t);\ldots;\vec{y}^M(t)]$ output of the network. The
network executes mapping $\vec{z}\mapsto L(\vec{z},\Theta)$ with network parameter $\Theta$.

\subsection{Neural Network Candidates ($L$)}
Choosing an RNN with feedforward ($\vec{F}$) and recurrent ($\vec{R}$) connections then the network assumes the form
\begin{equation}
\dot{\vec{y}}(\tau)=-\vec{y}(\tau)+\vec{F}\vec{z}(t)-\vec{R}\vec{y}(\tau)
\end{equation}
and thus, upon relaxation it solves the
\begin{equation}
\vec{y}(t)=(\vec{I}_D+\vec{R})^{-1}\vec{F}\vec{z}(t)=L(\vec{z}(t);\vec{F},\vec{R})\label{eq:network I-O2}
\end{equation}
input-output mapping \citep{base06blind,amari95recurrent}. Another natural choice is a network with feedforward
connections $\vec{W}$ that executes mapping
\begin{equation}
\vec{y}(t)=\vec{W}\vec{z}(t)=L(\vec{z}(t);\vec{W}).\label{eq:feedforward-network:IO}
\end{equation}

\subsection{Cost Function of \KISA}
The neural network estimates hidden sources $\vec{s}^m$ by non-linear ($\vec{f}$) decorrelation of $\vec{y}^m$s,
components of network output $\vec{y}$. Formally:

Let us denote the empirical $\vec{f}$-covariance matrix of $\vec{y}(t)$ and $\vec{y}^m(t)$ for function
$\vec{f}=[\vec{f}^1;\ldots;\vec{f}^M]$ over $[1,T]$ by
\begin{align}
     \vec{\Sigma}_{\K}(\vec{f},T)&=\widehat{cov}\left(\vec{f}\left[\pv_{\K}(\vec{y})\right],\vec{f}\left[\pv_{\K}(\vec{y})\right]\right),\label{eq:cov_C0}
     \\
     \vec{\Sigma}_{\K}^{i,j}(\vec{f},T)&=\widehat{cov}\left(\vec{f}^i\left[\pv_{\K}\left(\vec{y}^i\right)\right],\vec{f}^j\left[\pv_{\K}\left(\vec{y}^j\right)\right]\right),\label{eq:cov_C}
\end{align}
respectively, where $i,j=1,\ldots,M$, $\pv_{\R}(\vec{v})=\vec{v}$, $\pv_{\C}$ is the mapping
\begin{equation}
     \pv_{\C}:\C^L\ni\vec{v}\mapsto\vec{v}\otimes\left[\begin{array}{c}\Re(\cdot)\\\Im(\cdot)\end{array}\right]\in\R^{2L}.
\end{equation}
Here, $\Re(\cdot)$ [$\Im(\cdot)$] denotes the real (imaginary) part, $\otimes$ is the Kronecker-product. Then
minimization of the following non-negative cost function (in $\Theta$)
\begin{equation}
    Q_{\Theta}(\vec{f},T):=-\frac{1}{2}\log\left\{\frac{\det[\vec{\Sigma}_{\K}(\vec{f},T)]}{\prod_{m=1}^M\det[\vec{\Sigma}_{\K}^{m,m}(\vec{f},T)]}\right\}\label{IPA-cost}
\end{equation}
gives rise to \emph{pairwise}\footnote{We note that -- unlike in the the \mbox{1-dimensional} case, i.e., unlike for $d=1$ -- pairwise
independence is \emph{not }equivalent to mutual independence. Nonetheless, according to our numerical experiences it is an efficient
approximation.} \mbox{$\vec{f}$-uncorrelatedness}:
\begin{theorem}\label{thm:equiv}
For the separation carried out by the network minimizing cost function \eqref{IPA-cost}, the following statements are equivalent:
\begin{enumerate}
    \item[i)]\label{thm:equiv:i)}
        $\vec{f}$-uncorrelatedness: $\vec{\Sigma}_{\K}^{i,j}(\vec{f},T)=0$\quad($\forall i\ne j$).
    \item[ii)]
        $Q_{\Theta}$ is minimal: $Q_{\Theta}(\vec{f},T)=0$.
\end{enumerate}
\end{theorem}

\begin{proof2}[sketch]
The statement follows from the inequality related to the multi-dimensional Shannon differential entropy $H$: Let
$\vec{u}=[\vec{u}^1;\ldots;\vec{u}^M]\in\R^D$ $(\vec{u}^m\in\R^d)$ denote a random variable. Then
\begin{equation}
H(\vec{u}^1,\ldots,\vec{u}^M)\le\sum_{m=1}^MH(\vec{u}^m),
\end{equation}
and equality holds iff $\vec{u}^m$s are independent. Hint: one can choose $\vec{u}$ as a normal random variable with covariance
$\vec{\Sigma}_{\K}(\vec{f},T)$ and insert the expression of the entropy of normal variables.
\end{proof2}

\begin{note}
For the special case $\K=\R$, \mbox{$\Theta=(\vec{F},\vec{R})$}, \mbox{$\vec{f}(\vec{z})=\vec{z}$} and $d=1$, see \citep{base06blind}.
\end{note}

\begin{note}
Cost function $Q_{\Theta}$ of \eqref{IPA-cost} is attractive from the point of view of computing its gradient. This
gradient for the case of an RNN architecture [see Eq.~\eqref{eq:network I-O2}] may give rise to self-organization
\citep{base06blind}.
\end{note}

\begin{note}
For real random variables, the separation, which is aimed by cost function \eqref{IPA-cost}, can be related to the more general principle,
the Kernel Generalized Variance (KGV) technique \citep{bach02kernel}. This technique aims to separate the $\vec{y}^m$ components of
$\vec{y}$, the transformed form of input $\vec{z}$. To this end, KGV estimates mutual information $I(\vec{y}^1,\ldots,\vec{y}^M)$ in
Gaussian approximation\footnote{A complex variable is normal if its image using map $\pv_{\C}$ is real multivariate normal
\citep{eriksson06complex}. Thus, relation $I(\vec{y}^1,\ldots,\vec{y}^M)=I[\pv_{\C}(\vec{y}^1),\ldots,\pv_{\C}(\vec{y}^M)]\quad
\vec{y}^m\in\C^d$ extends the KGV based interpretation to the complex case, too [see Eqs.~\eqref{eq:cov_C0}-\eqref{eq:cov_C}].} by means of
the covariance matrix of variable $\vec{y}$.  Here, the transformation of the KGV technique is realized by the neural network parameterized
with variable $\Theta$ and by the function $\vec{f}$.
\end{note}

\begin{note}\label{note:KC}
We note that KGV is related to the kernel covariance (KC) method \citep{gretton03kernel}, which makes use of the supremum of
\mbox{1-dimensional} covariances as a measure of independence. Our approximation may also be improved by minimizing $Q_{\Theta}(\vec{f},T)$
on $\F(\ni\vec{f})$, i.e., on a set of functions.
\end{note}

\subsection{The \KISA Algorithm}\label{subsec:KISA-alg}
Below, our proposed \KISA method is introduced. A decomposition principle called \KISA Separation Theorem has been
formulated in \citep{szabo06separation}. It says that (under certain conditions) the \KISA task can be solved in 2
steps: In the first step, 1-dimensional \KICA estimation is executed that provides separation matrix
$\vec{W}_{\text{\KICA}}$ and estimated sources $\hat{\vec{s}}_{\text{\KICA}}$. In the second step, optimal permutation
of the \KICA elements ($\hat{\vec{s}}_{\text{\KICA}}$) is searched for, the \KICA elements are grouped.

This principle is adapted to linear feedforward neural networks [see, Eq.~\eqref{eq:feedforward-network:IO}]
here.\footnote{For the sake of simplicity we assume that all components have the same dimension, i.e., $d=d_m(\forall
m)$.} Separation matrix \mbox{$\vec{W}=\vec{W}_{\text{\KISA}}$} is searched in the form
\begin{equation}
    \vec{W}_{\text{\KISA}}=\vec{P}\vec{W}_{\text{\KICA}},
\end{equation}
where matrix $\vec{P}\in\R^{D\times D}$ denotes the desired permutation matrix. We search for the hidden sources
$\vec{s}^m$ by pairwise decorrelation of the components $\vec{y}^m$ of the output of the network using function
manifold $\F$ ($\F$: see, Note~\ref{note:KC}). Thus, given Theorem~\ref{thm:equiv}, our cost function is:
\begin{equation}
    Q(\F,T,\vec{P}):=\sum_{\vec{f}\in\F}\left\|\vec{M}\circ\vec{\Sigma}_{\K}(\vec{f},T,\vec{P})\right\|^2\rightarrow\min_{\vec{P}}.\\
\end{equation}
Here: (i) $\F$ denotes a set of functions, each function $\R^D\mapsto \R^D$ (if $\K=\C$ then $\R^{2D}\mapsto \R^{2D}$), and each function
acts on each coordinate separately, (ii) $\circ$ denotes pointwise multiplication (Hadamard product), (iii) $\vec{M}$ masks according to
the subspaces \{\mbox{$\vec{M}=\vec{E}_D-\vec{I}_M\otimes \vec{E}_d$}, where all elements of matrix \mbox{$\vec{E}_D\in\R^{D\times D}$} and
\mbox{$\vec{E}_d\in\R^{d\times d}$} are equal to 1 [if $\K=\C$ then  $\vec{E}_D$ ($\vec{E}_d$) is replaced by $\vec{E}_{2D}$
($\vec{E}_{2d}$)]\}, (iv) $\left\|\cdot\right\|^2$ denotes the square of the Frobenius norm (sum of squares of the elements), (v) in
$\vec{\Sigma}_K(\vec{f},T,\vec{P}$), \mbox{$\vec{y}=\vec{P}\hat{\vec{s}}_{\text{\KICA}}$}, and (vi) $\vec{P}$ is the $D\times D$
permutation matrix to be determined.

Greedy permutation search is applied: 2 coordinates of different subspace are exchanged if this change lowers cost function
$Q(\F,T,\cdot)$. Note: Greedy search could be replaced by a \emph{global} one for a higher computational burden \citep{szabo06cross}.
Table~\ref{tab:KISA-pseudocode} contains the pseudocode of our technique.

\begin{table}
  \centering
  \caption{\KISA Algorithm - pseudocode}\label{tab:KISA-pseudocode}
  \begin{minipage}{\textwidth}
  \begin{tabular}{|l|}
        \hline
        \textbf{Input of the algorithm}\\
        \verb|  |observation: $\{\mathbf{z}(t)\}_{t=1,\ldots,T}$\\
        \textbf{Optimization}\footnote{Let $\G^1, \ldots, \G^M$ denote the indices of the $1^{st}, \ldots , M^{th}$\\
             subspaces, i.e., $\G^m:=\{(m-1)d+1,\ldots,md\}$, and\\
             permutation matrix $\vec{P}_{pq}$ exchanges coordinates $p$ and $q$.}\\
        \verb|  |\textbf{\KICA}: on whitened observation $\vec{z}$, \\
        \verb|        |$\Rightarrow$ $\hat{\vec{s}}_{\text{\KICA}}$ estimation\\
        \verb|  |\textbf{Permutation search}\\
        \verb|    |$\vec{P}:=\vec{I}_D$\\
        \verb|    |repeat\\
        \verb|      |sequentially for $\forall p\in\G^{m_1},q\in\G^{m_2}$\\
        \verb|                     |$(m_1\ne m_2):$\\
        \verb|        |if $Q(\F,T,\vec{P}_{pq}\vec{P})<Q(\F,T,\vec{P})$\\
        \verb|          |$\vec{P}:=\vec{P}_{pq}\vec{P}$\\
        \verb|        |end\\
        \verb|    |until $Q(\F,T,\cdot)$ decreases in the \emph{sweep} above\\
        \textbf{Estimation}\\
        \verb|  |$\hat{\vec{s}}_{\text{\KISA}}=\vec{P}\hat{\vec{s}}_{\text{\KICA}}$\\
        \hline
  \end{tabular}
  \end{minipage}
\end{table}

\section{Illustrations}\label{sec:illustrations}
The \KISA identification algorithm of Section~\ref{subsec:KISA-alg} is illustrated below (due to the lack of space,
illustrations are provided for $\K=\R$ only). Test cases are introduced in Section~\ref{subsec:databases}. The quality
of the solutions will be measured by the normalized Amari-distance (Section~\ref{subsec:amaridist}). Numerical results
are provided in Section~\ref{subsec:simulations}.

\subsection{Databases}\label{subsec:databases}
Two databases were defined to study our algorithm. The databases are illustrated in Fig.~\ref{fig:database-Aw} and
Fig.~\ref{fig:database-d-spherical}.

\subsubsection{The $A\omega$ Database}
\begin{figure}%
\centering%
\includegraphics[width=3cm]{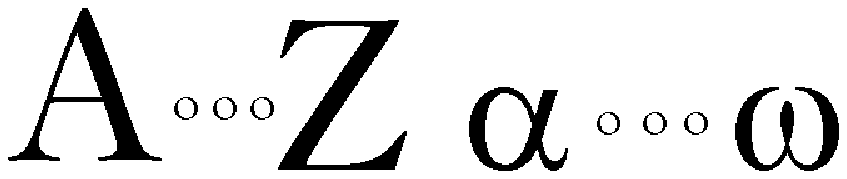}
\caption[]{Database $A\omega$. 100-dimensional task of 50 pieces of 2-dimensional components ($D=100$, $M=50$, \mbox{$d=2$}).
Hidden sources are uniformly distributed variables on the letters of the English and the Greek alphabets.}%
\label{fig:database-Aw}%
\end{figure}

Here, hidden sources $\vec{s}^m$ are uniform distributions defined by the 2-dimensional images ($d=2$) of letters on
A-Z and $\alpha-\omega$. This is called database $A\omega$, which has 50 components ($M=50$), see
Fig.~\ref{fig:database-Aw} for an illustration. This test falls outside of the (known) validity domain of the \RISA
Separation Theorem.

\subsubsection{The $d$-Spherical Database}
\begin{figure}%
\centering%
\subfloat[][]{\includegraphics[trim=100 27 70 40,scale=0.17]{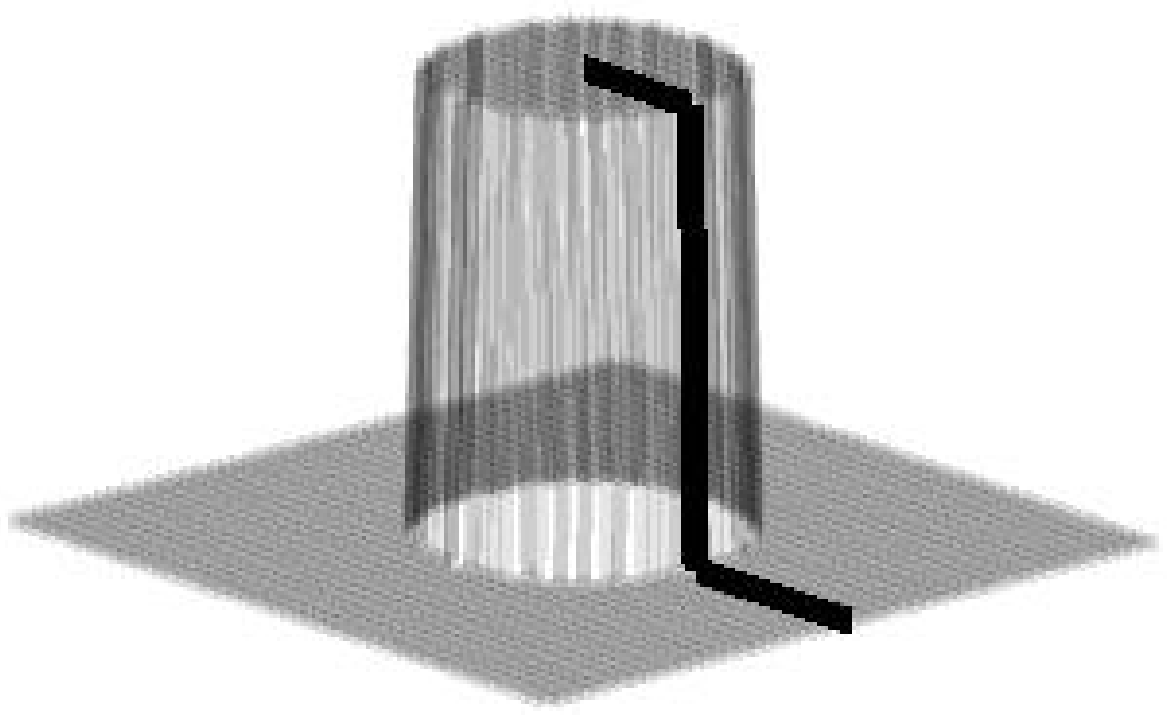}}\hfill%
\subfloat[][]{\includegraphics[trim=85 30 67 40,scale=0.15]{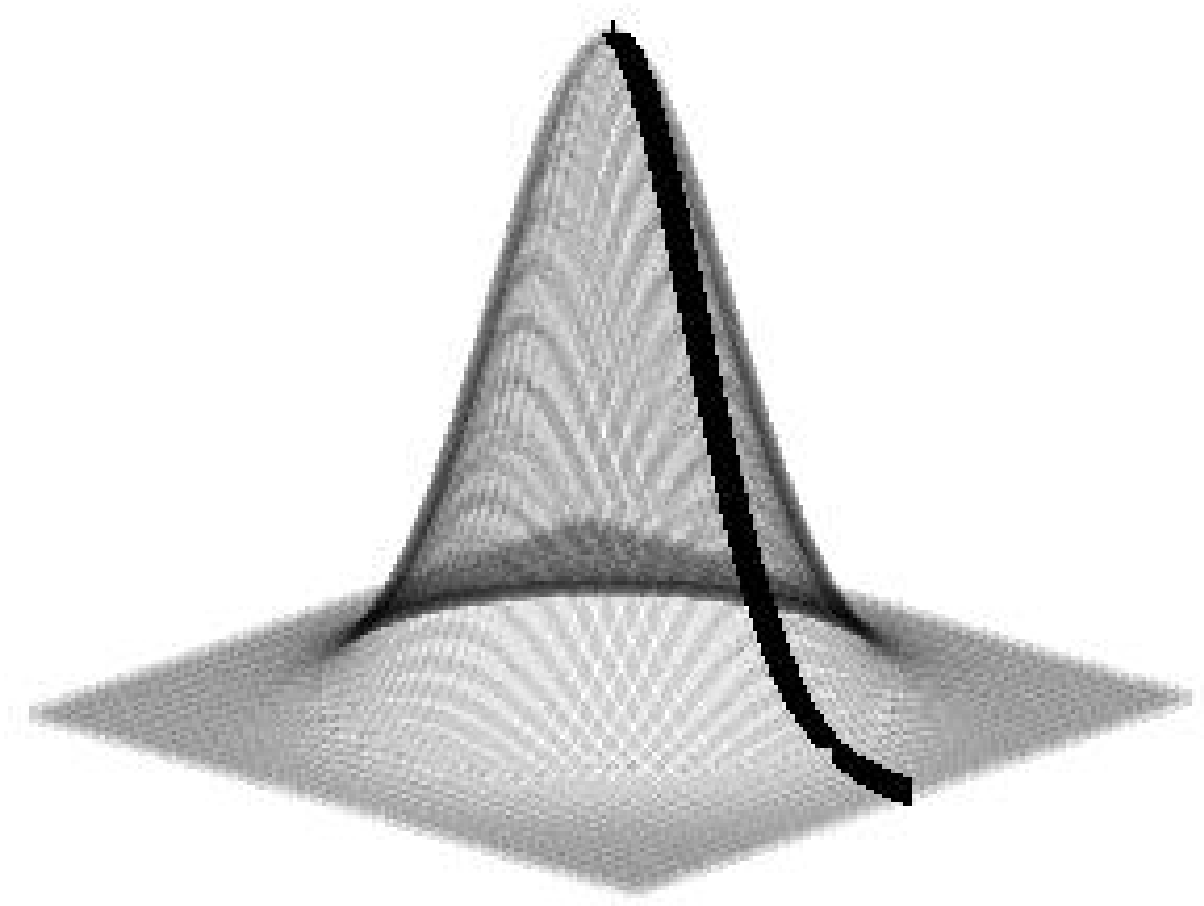}}\hfill%
\subfloat[][]{\includegraphics[trim=90 30 75 40,scale=0.16]{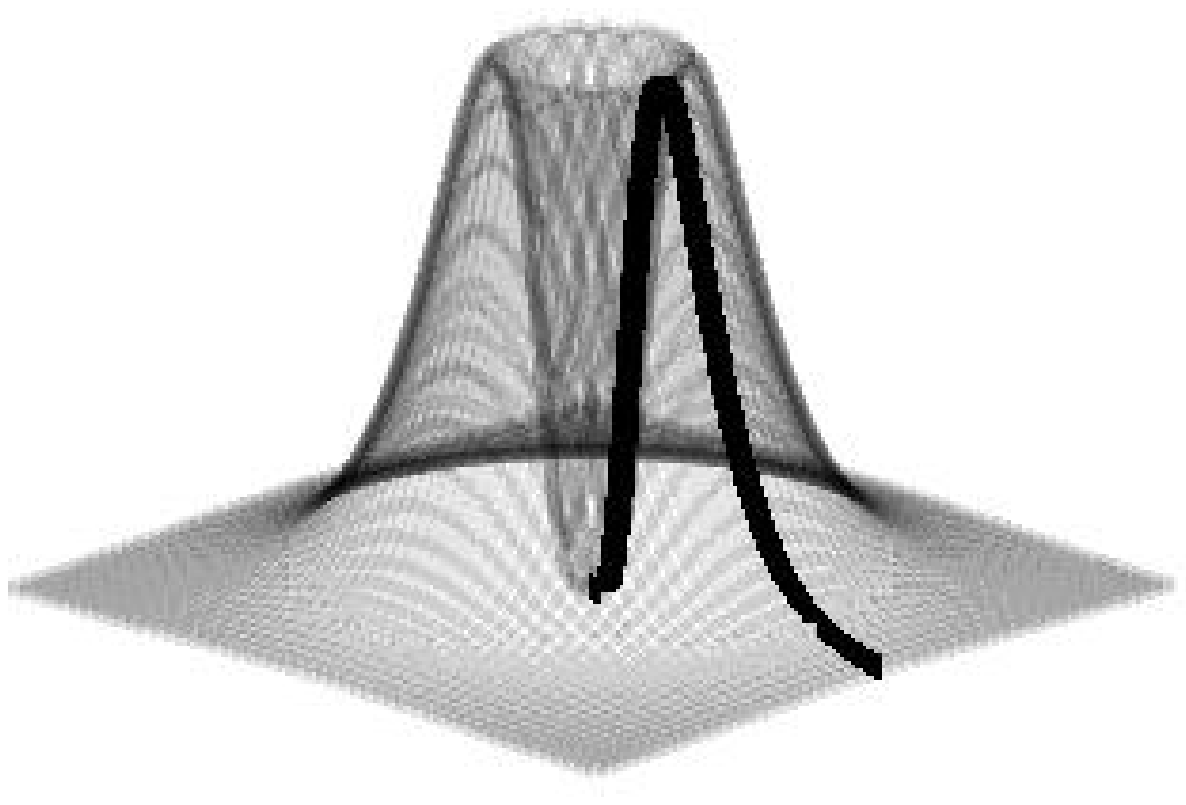}}\hfill%
\caption[]{Database $d$-spherical. Stochastic representation of the 3 ($M=3$) hidden sources. (a): $\rho$ is
uniform on $[0,1]$, (b): $\rho$ is exponential with parameter $\mu=1$, and (c): $\rho$ is lognormal with parameters $\mu=0, \sigma=1$, respectively.}%
\label{fig:database-d-spherical}%
\end{figure}

Here, hidden sources $\vec{s}^m$ are spherically symmetric random variables that have representation of the form
\mbox{$\vec{v}\stackrel{\mathrm{distr}}{=}\rho\vec{u}^{(d)}$}, where $\vec{u}^{(d)}$ is uniformly distributed on the $d$-dimensional unit
sphere, and $\rho$ is a non-negative scalar random variable independent of $\vec{u}^{(d)}$ ($\stackrel{\mathrm{distr}}{=}$ denotes equality
in distribution). This \emph{$d$-spherical} database: (i) can be scaled in dimension $d$, (ii) satisfies conditions of the \RISA Separation
Theorem, and (iii) can be defined by $\rho$. (See \citep{fang90symmetric,frahm04generalized} for spherical variables.) Our choices for
$\rho$ are shown in Fig.~\ref{fig:database-d-spherical}.

\subsection{Normalized Amari-distance}\label{subsec:amaridist}
The precision of our algorithm was measured by the normalized Amari-distance as follows. The optimal estimation of the \RISA model provides
matrix \mbox{$\vec{B}:=\vec{W}\vec{A}\in\R^{D\times D}$}, a block-permutation matrix made of $d\times d$ sized blocks. Let us decompose
matrix $\vec{B}\in\R^{D\times D}$ into $d\times d$ blocks: \mbox{$\vec{B}=\left[\vec{B}^{i,j}\right]_{i,j=1,\ldots,M}$}. Let $b^{i,j}$
denote the sum of the absolute values of the elements of matrix \mbox{$\vec{B}^{i,j}\in\R^{d\times d}$}. Then the \emph{normalized} version
of the Amari-distance \citep{amari96new} (as it was introduced in \citep{szabo06cross} for \RISA) is defined as:
\begin{eqnarray}
    r(\vec{B})&:=&\frac{1}{2M(M-1)}\left[\sum_{i=1}^M\left(\frac{
    \sum_{j=1}^Mb^{ij}}{\max_jb^{ij}}-1\right)\right. + \nonumber \\
    &&\qquad \qquad\quad\left.\sum_{j=1}^M\left(\frac{ \sum_{i=1}^Mb^{ij}}{\max_ib^{ij}}-1\right)\right].
\end{eqnarray}
For matrix $\vec{B}$ we have that $0\le r(\vec{B})\le 1$, and $r(\vec{B})=0$ if, and only if $\vec{B}$ is a block-permutation matrix with
$d\times d$ sized blocks. Thus, $r=0$ corresponds to perfect estimation (0\% error), $r=1$ is the worst estimation (100\% error). This
performance measure can be used for $\K=\C$, too.

\subsection{Simulations}\label{subsec:simulations}
Results on databases \emph{$A\omega$} and \emph{$d$-spherical} are provided here.  In our simulations, sample number of observations
$\vec{z}(t)$ changed: $1000 \le T \le 30000$. Mixing matrix $\vec{A}$ was chosen randomly from the orthogonal group. Manifold $\F$ was
$\F:=\{\vec{z}\mapsto \cos(\vec{z}), \vec{z}\mapsto \cos(2\vec{z})\}$ (functions operated on coordinates separately). Scaling properties of
the approximation were studied for database \mbox{$d$-spherical} by changing the value of $d$ between $20$ and $110$ [i.e., the number of
subspaces ($M$) was fixed, but the dimension of the subspaces was increased.] For each parameters [$T$ for database $A\omega$, $(T,d)$ for
database $d$-spherical] ten experiments were averaged. Qualities of the solutions were measured by the Amari-error (see
Section~\ref{subsec:amaridist}). We have chosen FastICA \citep{hyvarinen97fast} for the \RICA module (see Table~\ref{tab:KISA-pseudocode}).

\begin{figure}%
\centering%
\subfloat[][]{\includegraphics[width=6.1cm]{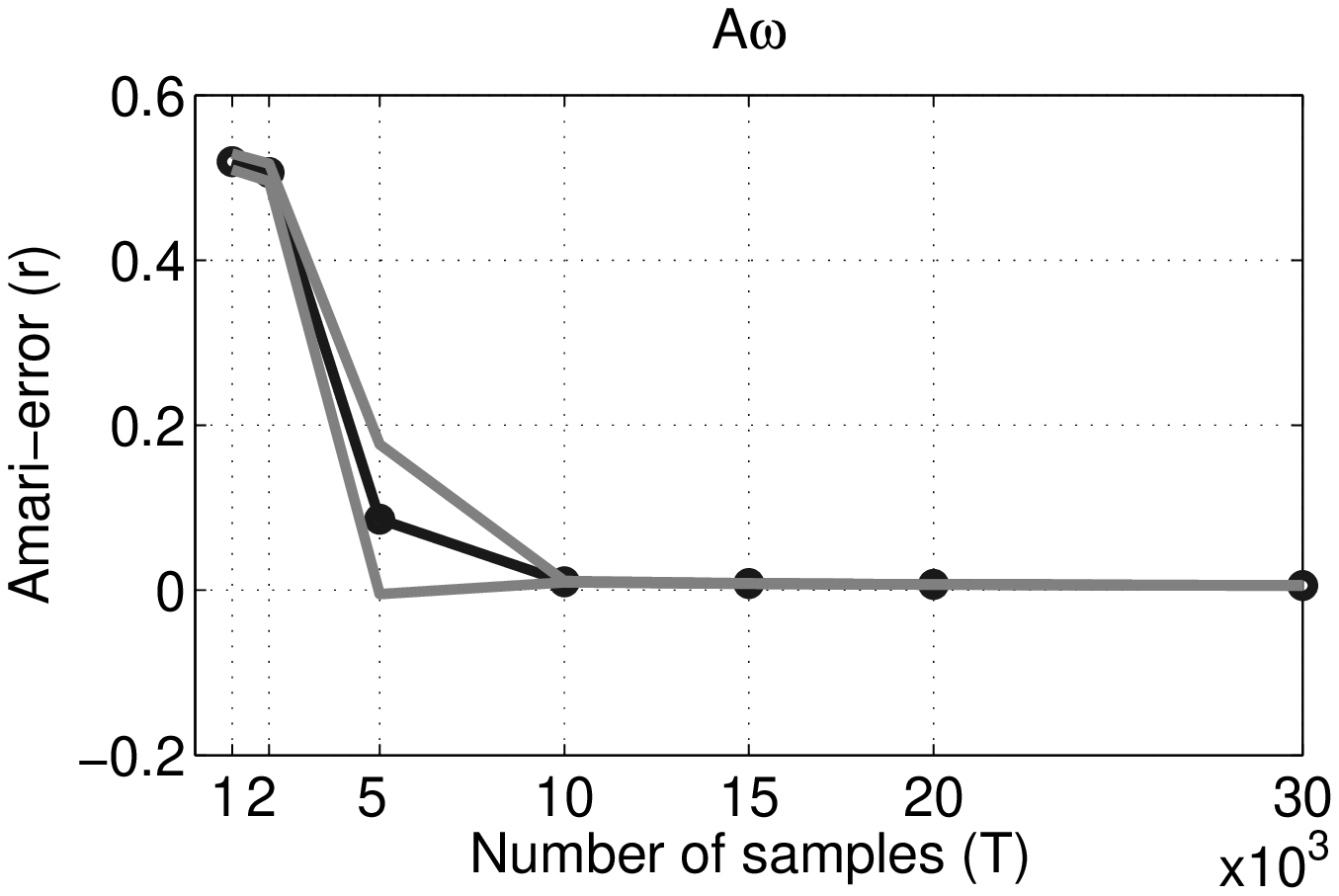}}%
\subfloat[][]{\includegraphics[width=2.3cm]{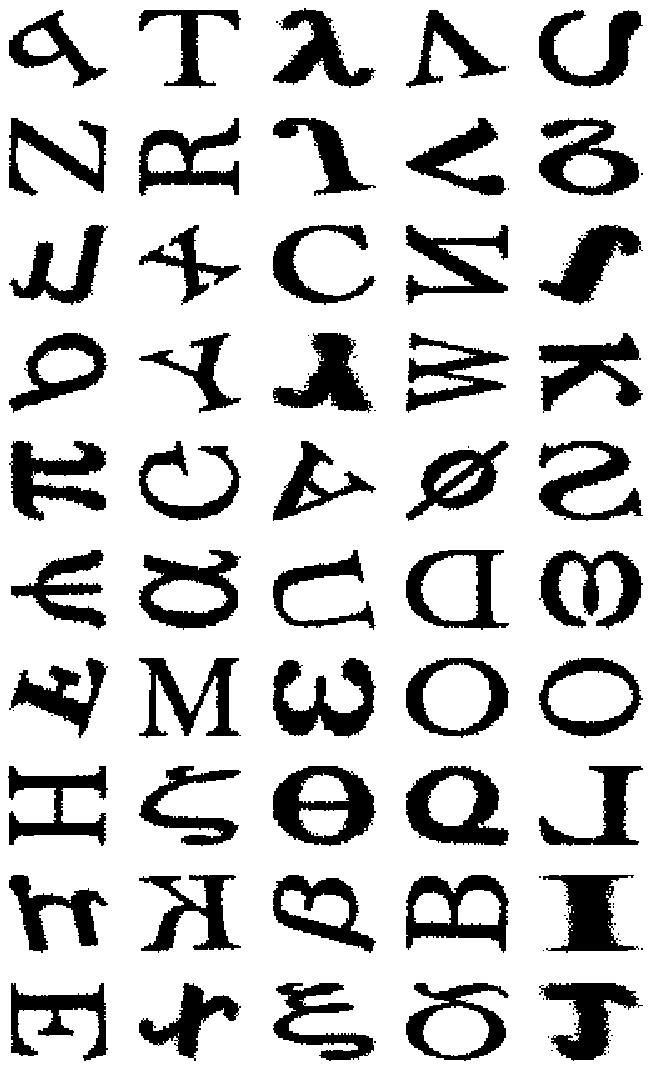}}\\%
\caption[]{Estimations on database $A\omega$. (a) Amari-error as a function of the number of samples.
Average$\pm$deviation for $30000$ samples: $0.58\%\pm0.04$, (b) estimation with average error for $30000$ samples:
the hidden components are recovered up to permutation and orthogonal transformation (\RISA ambiguity).}%
\label{fig:demo:A-omega}%
\end{figure}

\begin{figure}
\centering
\includegraphics[width=6.1cm]{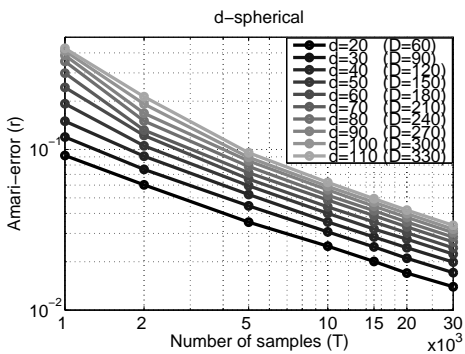}
\caption{Estimations of database $d$-spherical: Amari-error as a function of the number of samples on loglog scale for
different dimensional ($d$) subspaces. Task dimension: $D$. Errors are approximately linear, so they scale according to
power law, like $r(T)\propto T^{-c}$ $(c>0)$. For numerical values, see Table~\ref{tab:amari-dists-d-spherical}.}
\label{fig:amari-vs-T:d-spherical}
\end{figure}

Precision of our method is shown: (i) for database $A\omega$ in Fig.~\ref{fig:demo:A-omega} as a function of sample
number, (ii) for database $d$-spherical in Fig.~\ref{fig:amari-vs-T:d-spherical} as a function of sample number and
source dimension ($d$) (for details, see Table~\ref{fig:amari-vs-T:d-spherical}). The figures demonstrate that the
algorithm was able to uncover the hidden components with high precision. In the case of database $d$-spherical the
Amari error decreases according to power law $r(T)\propto T^{-c}$ $(c>0)$.

In our numerical simulations, the number of sweeps before the iteration of the permutation optimization stopped (see
Table~\ref{tab:KISA-pseudocode}) varied between 2 and 6.

\begin{table}
    \centering
    \caption{Amari-error for database $d$-spherical, for different $d$ values: average $\pm$ deviation. Number of samples: $T=30000$.}
    \label{tab:amari-dists-d-spherical}
    \begin{tabular}{|c|c|c|}
    \hline
        $d=20$ & $d=30$ &$d=40$\\
    \hline
        $1.40\%$ $(\pm 0.03)$ & $1.71\%$ $(\pm 0.03)$& $1.99\%$ $(\pm 0.03)$\\
    \hline\hline
       $d=50$ & $d=60$& $d=70$\\
    \hline
          $2.23\%$ $(\pm 0.03)$& $2.44\%$ $(\pm 0.03)$ & $2.65\%$ $(\pm 0.03)$\\
    \hline\hline
        $d=80$& $d=90$& $d=100$\\
    \hline
         $2.85\%$ $(\pm 0.03)$& $3.03\%$ $(\pm 0.04)$&$3.19\%$ $(\pm 0.02)$\\
    \hline\hline
        $d=110$&\multicolumn{2}{c}{}\\
    \cline{1-1}
        $3.37\%$ $(\pm 0.03)$&\multicolumn{2}{c}{}\\
    \cline{1-1}
    \end{tabular}
\end{table}


\begin{thebibliography}{13}
\providecommand{\natexlab}[1]{#1} \providecommand{\url}[1]{\texttt{#1}} \expandafter\ifx\csname urlstyle\endcsname\relax
  \providecommand{\doi}[1]{doi: #1}\else
  \providecommand{\doi}{doi: \begingroup \urlstyle{rm}\Url}\fi

\bibitem[Amari et~al.(1995)Amari, Cichocki, and Yang]{amari95recurrent}
S.~Amari, A.~Cichocki, and H.H. Yang.
\newblock Recurrent neural networks for blind separation of sources.
\newblock NOLTA'95, pp. 37-42, 1995.

\bibitem[Amari et~al.(1996)Amari, Cichocki, and Yang]{amari96new}
S.~Amari, A.~Cichocki, and H.~Yang.
\newblock A new learning algorithm for blind signal separation.
\newblock \emph{Advances in Neural Information Processing Systems}, 1996.

\bibitem[Bach and Jordan(2002)]{bach02kernel}
Francis~R. Bach and Michael~I. Jordan.
\newblock Kernel {ICA}.
\newblock \emph{JMLR}, 3:\penalty0 1--48, 2002.

\bibitem[Eriksson and Koivunen(2006)]{eriksson06complex}
Jan Eriksson and Visa Koivunen.
\newblock Complex random vectors and {ICA} models: Identifiability, uniqueness
  and separability.
\newblock \emph{{IEEE} {TIT}}, 52\penalty0 (3), 2006.

\bibitem[Fang et~al.(1990)Fang, Kotz, and Ng]{fang90symmetric}
Kai-Tai Fang, Samuel Kotz, and Kai~Wang Ng.
\newblock \emph{Symmetric multivariate and related distributions}.
\newblock Chapman and Hall, 1990.

\bibitem[Frahm(2004)]{frahm04generalized}
Gabriel Frahm.
\newblock \emph{Generalized elliptical distributions: Theory and applications}.
\newblock PhD thesis, University of Köln, 2004.

\bibitem[Gretton et~al.(2003)Gretton, Herbrich, and Smola]{gretton03kernel}
A.~Gretton, R.~Herbrich, and A.~Smola.
\newblock The kernel mutual information.
\newblock In \emph{IEEE ICASSP}, volume~4, pages 880--883, 2003.

\bibitem[Hyv{\"a}rinen and Oja(1997)]{hyvarinen97fast}
A.~Hyv{\"a}rinen and E.~Oja.
\newblock A fast fixed-point algorithm for independent component analysis.
\newblock \emph{Neural Computation}, 9\penalty0 (7):\penalty0 1483--1492, 1997.

\bibitem[Meyer-B{\"a}se et~al.(2006)Meyer-B{\"a}se, Gruber, Theis, and
  Foo]{base06blind}
Anke Meyer-B{\"a}se, Peter Gruber, Fabian Theis, and Simon Foo.
\newblock Blind source separation based on self-organizing neural network.
\newblock \emph{Engng. Appl. of AI}, 19:\penalty0 305--311, 2006.

\bibitem[P\'oczos and L\H{o}rincz(2006)]{poczos06noncombinatorial}
B.~P\'oczos and A.~L\H{o}rincz.
\newblock Non-combinatorial estimation of independent {AR} sources.
\newblock \emph{Neurocomp.}, 2006.
\newblock accepted.

\bibitem[Szab\'o et~al.(2006{\natexlab{a}})Szab\'o, P\'oczos, and
  L\H{o}rincz]{szabo06cross}
Z.~Szab\'o, B.~P\'oczos, and A.~L\H{o}rincz.
\newblock Cross-entropy optimization for independent process analysis.
\newblock In \emph{Proc. of {ICA}}, LNCS 3889, pages 909--916. Springer-Verlag,
  2006{\natexlab{a}}.

\bibitem[Szab\'o et~al.(2006{\natexlab{b}})Szab\'o, P\'oczos, and L{\H
  o}rincz]{szabo06separation}
Z.~Szab\'o, B.~P\'oczos, and A.~L{\H o}rincz.
\newblock Separation theorem for $\mathbb{K}$-independent subspace analysis
  with sufficient conditions.
\newblock Technical report, E\"otv\"os Lor\'and University, Budapest,
  2006{\natexlab{b}}.
\newblock http://arxiv.org/abs/math.ST/0608100.

\bibitem[Theis(2004)]{theis04uniqueness1}
F.J. Theis.
\newblock Uniqueness of complex and multidimensional independent component
  analysis.
\newblock \emph{Signal Processing}, 84\penalty0 (5):\penalty0 951--956, 2004.
\end{thebibliography}
\end{document}